\documentclass{amsart}

\usepackage{amssymb, amsmath,amsthm}
\usepackage{amsaddr}
\usepackage{graphicx}
\usepackage{mathrsfs}
\usepackage{color}
\usepackage{esint}
\usepackage{enumerate}
\usepackage[colorlinks=true, linkcolor=blue,citecolor=blue]{hyperref}

\newtheorem{theorem}{Theorem}[section]

\newtheorem{algorithm}{Algorithm}

\newtheorem{lemma}[theorem]{Lemma}

\theoremstyle{definition}
\newtheorem{definition}[theorem]{Definition}
\newtheorem{remark}{Remark}

\newtheorem{assumption}{Assumption}
\numberwithin{equation}{section}

\title[Inverse Problem for the Non-self-adjoint Schr\"odinger Equation]{Inverse Problem for the Schr\"odinger Equation with Non-self-adjoint Matrix Potential}
\author[S A Avdonin]{Sergei Avdonin}
\address{Department of Mathematics and Statistics, University of Alaska Fairbanks, AK 99775-6660, USA}
\author[A S Mikhaylov]{Alexander Mikhaylov}
\address{St. Petersburg   Department   of   V.A. Steklov Institute   of   Mathematics
of   the   Russian   Academy   of   Sciences, 27, Fontanka,
191023, St. Petersburg, Russia and Saint Petersburg State
University, St.Petersburg State University, 7/9 Universitetskaya
nab., St. Petersburg, 199034, Russia.}
\author[V S Mikhaylov]{Victor Mikhaylov}
\address{St. Petersburg   Department   of   V.A. Steklov    Institute   of   Mathematics
of   the   Russian Academy   of   Sciences, 27, Fontanka, 191023,
St. Petersburg, Russia}
\author[J Park]{Jeff Park}
\address{Department of Mathematics and Statistics, University of Alaska Fairbanks, AK 99775-6660, USA}

\email{saavdonin@alaska.edu}
\email{mikhaylov@pdmi.ras.ru}
\email{vsmikhaylov@pdmi.ras.ru}
\email{jcpark2@alaska.edu}

\begin{document}

\begin{abstract}
We consider the dynamical system with boundary control for the
vector Schr\"odinger equation on the interval with a
non-self-adjoint matrix potential. For this system, we study the
inverse problem of recovering the matrix potential from the
dynamical Dirichlet--to--Neumann operator. We first provide a
method to recover spectral data for an abstract system from
dynamic data and apply it to the Schr\"odinger equation. We then
develop a strategy for solving the inverse problem for the
Schr\"odinger equation using this method with other techniques of
the Boundary control method.
\end{abstract}

\maketitle


\section{Introduction}

For this work, we consider the following dynamical system
\begin{equation} \label{Dsyst}
\left\{ \begin{array}{ll}
iu_t - u_{xx} + Q(x)u = 0, &\quad 0 \leq x \leq \ell, \: 0 < t < T,\\
u(x,0) = u_t (x,0) = 0, &\quad 0 \leq x \leq \ell,\\
u(0,t) = f(t), \: u(\ell,t) = 0, &\quad 0 < t < T,
\end{array} \right.
\end{equation}
where $\ell > 0$, $T>0$ are given, $Q \in C^2 \big( (0,\ell);
\mathbb{R}^{N \times N} \big)$, $Q \neq Q^*$, is a matrix
potential. The vector function $f \in L^2 \big( (0,T);
\mathbb{R}^N \big)$ is referred to as the \emph{boundary control}.
The solution to (\ref{Dsyst}) is denoted $u^f$. We introduce the
response operator $R^T$ by
\begin{align*}
R^T :\: &L^2 \big( (0,T);\mathbb{R}^N \big) \rightarrow L^2 \big( (0,T);\mathbb{R}^N \big),\\
&(R^T f) (t) = u_x^f (0,t), \quad 0 < t < T. \notag
\end{align*}
The inverse problem is to recover $Q(x)$ for $0 < x < \ell$ from $R^T$.

The most common method to solve the inverse problem involves
recovering spectral data---eigenvalues and traces of
eigenfunctions---from the dynamical data, $R^T$, and then solving
the resulting spectral problem (see \cite{alp2002}). The
connections among the different types of data (dynamical,
spectral, scattering) are an important topic in the theory of
inverse problems, see \cite{b2003, b2001, kklm2004,mm2017} to
mention a few. When the system is spectrally controllable, the
variational method is used to obtain spectral data (see
\cite{alp2002, b2003} for details). This approach, which is used
in \cite{b2003, b1995, abr1997}, is based on the Boundary control
(BC) method and relies on the operator being self--adjoint.

In the present paper, we are considering non-self-adjoint
operators and the variational method is not applicable. However, we will
still follow a similar strategy, i.e., recovering spectral
data from the dynamical data. Instead of using the variational
approach, we will use a different method, proposed in
\cite{agm2010}, that also uses spectral controllability of the
underlying dynamical system.

We begin by considering the dynamical control system in an
abstract setting. Recently in \cite{mm2012, amm2015}, the authors
studied the same questions for the adjoint system, a dynamical
system with observation. They derived the equations of the
boundary control (BC) method for this system (see \cite{b2007} for
explanation of these equations). Using these equations, they
treated the one-dimensional inverse source problem and the
dynamical inverse problem with one measurement for the
Schr\"odinger equation in \cite{am2011, amr2014}. In the case of a
system with observation, only one measurement is available,
however, we consider the response operator on an interval.

As mentioned, Section 2 will introduce the abstract system and
derive the equations of the BC method. Sections 3 and 4 are
concerned with recovering spectral data for an operator with a
simple spectrum. Section 5 solves the same problem for an operator
whose spectrum is not simple. In the last section, we prove
spectral controllability for the Schr\"odinger system
(\ref{Dsyst}) and recover the matrix potential $Q$.


\section{Equations of the BC method}

In this section, we consider an abstract dynamical control system.
Let $H$ and $Y$ be Hilbert spaces, and $A$ an operator in $H$ that
is not necessarily self-adjoint. We consider the dynamical system
in $H$:
\begin{equation} \label{ED}
\left\{ \begin{array}{ll}
iu_t - Au = Bf, &\quad t>0,\\
u(0) = 0,
\end{array} \right.
\end{equation}
where $B : Y \to H$ is an input operator. We define the observation operator $O$ by
\begin{align*}
O: &H \rightarrow Y\\
&O = B^*.
\end{align*}
We will fix $T>0$ and denote the solution to (\ref{ED}) by $u^f$ for $0 < t < T$.
We will define the response operator $R^T$ by
\begin{align*}
R^T :\: &L^2 \big( (0,T);Y\big) \rightarrow L^2 \big( 0,T);Y\big),\\
&\big( R^T f \big) (t) := \big( Ou^f \big)(t).
\end{align*}
Hence, $R^T$ is the output of the system.

Let $A^*$ denote the operator adjoint to $A$. Along with System
(\ref{ED}), we consider the following dynamical control system:
\begin{equation} \label{ED1}
\left\{ \begin{array}{ll} iv_t + A^*v = -Bg, &\quad t > 0,\\
v(0) = 0, \end{array} \right.
\end{equation}
and denote its solution by $v^g$. The response operator for this
system will be denoted $R_\#^T$, where $\big( R_\#^T g \big) (t)
:= \big( O v^g \big) (t)$, $t \in (0,T)$. For now, we will denote
$\mathcal{F}^T = L^2 \big( (0,T);Y\big)$. It is not difficult to
show the relationship between the response operators of Systems
(\ref{ED}) and (\ref{ED1}). We first introduce the operator $J^T$
in $L^2 \big( (0,T); Y \big)$ by the rule
\begin{equation}
\label{JT} \big( J^T f \big) (t) := f(T-t), \quad 0 \leq t \leq T.
\end{equation}
\begin{lemma} The following identity holds.
\begin{equation} \label{R_lnk}
\big( R_\#^T \big)^* J^T = J^T R^T.
\end{equation}
\begin{proof}
We introduce the function $w = v(T-t)$, which is a solution to
\[ \left\{ \begin{array}{ll}
iw_t - A^* w = B g(T-t), &\quad t >0,\\
w(T) = 0.
\end{array} \right. \]
Then we evaluate
\begin{align*}
\int_0^T \big( iu_t^f - Au^f, w^g \big)_H \: dt &= \int_0^T \big( Bf, w^g \big)_H \: dt\\
&= \int_0^T \big( f, Ow^g \big)_Y \: dt\\
&= \int_0^T \big( f, (Ov^g)(T-t) \big)_Y \: dt\\
&= \big( f(T-t), R_\#^T g \big)_{\mathcal{F}^T}\\
&= \Big( \big( R_\#^T\big)^* J^T f, g \Big)_{\mathcal{F}^T}.
\end{align*}
On the other hand, using integration by parts yields
\begin{align*}
\int_0^T \big(iu_t^f - Au^f, w^g \big)_H \: dt &= \int_0^T (u^f, iw_t^g - A^* w^g)_H \: dt\\
&= \int_0^T \big( u^f, Bg(T-t)\big)_H \: dt\\
&= \int_0^T \big( Ou^f, g(T-t))_Y \: dt\\
&= \big( R^T f, g(T-t) \big)_{\mathcal{F}^T}\\
&= \Big( J^T R^T f, g \Big)_{\mathcal{F}^T}.
\end{align*}
Comparing the last expressions with the fact that $f$ and $g$ are arbitrary completes the proof.
\end{proof}
\end{lemma}

For Systems (\ref{ED}) and (\ref{ED1}), we introduce the \textit{control operators}
\begin{align*}
W^T &: \mathcal{F}^T \to H, \quad W^T f := u^f (T),\\
W_\#^T &: \mathcal{F}^T \to H, \quad W_\#^T g := v^g (T).
\end{align*}
From the control operators, we introduce the \emph{connecting
operator} $C^T : \mathcal{F}^T \to \mathcal{F}^T$ by its quadratic
form
\[ (C^T f,g)_{\mathcal{F}^T} = \big( u^f (T), v^g (T) \big)_H = \big( W^T f, W_\#^T g \big)_H. \]

It is an important fact in the BC method that $C^T$ can be
expressed in terms of the inverse data. For this, we use the
operator $J^{2T}$ in $\mathcal{F}^{2T}$ defined in (\ref{JT}) and
$Z^T : \mathcal{F}^T \to \mathcal{F}^{2T}$ defined by the rule
\[ \big( Z^T f \big) (t) := \left\{ \begin{array}{ll} f(t), &0 \leq t \leq T,\\
0, & T < t \leq 2T.\end{array} \right.
\]
\begin{lemma} The following representation holds.
\begin{equation} \label{CT_repr}
C^T = -i \big(Z^T \big)^* J^{2T} R^{2T} Z^T.
\end{equation}
\begin{proof}
We introduce the Blagoveschenskii function defined by
\[ \psi (s,t) = \big( u^f (s), v^g (t) \big)_H \]
and evaluate
\begin{align*}
\psi_s (s,t) &= \big( -iAu^f (s) - iB f(s), v^g (t) \big)_H \\
&= \big( u^f (s), iA^* v^g (t) \big)_H - \big( if(s), Ov^g (t) \big)_Y,\\
\psi_t (s,t) &= \big( u^f (s), iA^*v^g (t) + iBg(t) \big)_H \\
&= \big( u^f (s), iA^* v^g (t) \big)_H + \big( Ou^f (s), i g(t) \big)_Y.
\end{align*}
Thus, $\psi(s,t)$ satisfies
\begin{align*}
\psi_t (s,t) - \psi_s (s,t) &= -i \bigg( \big(R^T f \big) (s), g(t) \bigg)_Y + i \bigg( f(s), \big(R_\#^T g \big) (t) \bigg)_Y =: h(s,t),\\
\psi(0,s)  &= 0.
\end{align*}
Integrating this equation yields
\[ \psi (s,t) = \int_0^t h(s+t-\eta, \eta) \: d\eta, \]
where we then set $f(t) = 0$ for $t \notin (0,T)$ and get
\begin{align*}
\big( C^T f,g \big)_{\mathcal{F}^T} = \psi (T,T) &= \int_0^T h(2T -\eta, \eta ) \: d\eta\\
&= -i \int_0^T \bigg( \big(R^{2T} f \big) (2T-\eta), g(\eta) \bigg)_Y \: d\eta,
\end{align*}
which completes the proof.
\end{proof}
\end{lemma}


\section{The Spectral Problem for the Simple Case -- algebraically simple spectrum}
In what follows, we will assume that $A$ satisfies the following:
\begin{assumption} \label{A1} \mbox{}
\begin{enumerate}[(a)]
\item The spectrum of $A$ is simple, i.e., it consists of
(infinitely many) eigenvalues with algebraic multiplicity one. We
denote them by $\{\lambda_k\}_{k=1}^\infty$ and the adjoint
operator $A^*$ has spectrum $\{
\overline{\lambda_k}\}_{k=1}^\infty$. \item The eigenfunctions of
$A$ form a Riesz basis in $H$, denoted $\{ \varphi_k
\}_{k=1}^\infty$, the basis of $A^*$ we denote by $\{
\psi_k\}_{k=1}^\infty$, and the property $\left(\varphi_k,
\psi_l\right)_H = \delta_{kl}$ holds. \item Systems (\ref{ED}) and
(\ref{ED1}) are spectrally controllable, i.e., there exist
controls $f_k, g_k \in H_0^1 \big( (0,T);Y \big)$ such that $W^T
f_k = \varphi_k$ and $W_\#^T g_k = \psi_k$.
\end{enumerate}
\end{assumption}

By dot, we denote differentiation with respect to $t$. We formulate the main result.
\begin{theorem} \label{thm1}
If $A$ satisfies Assumption \ref{A1}, then the spectrum of $A$ and
(non-normalized) controls $f_k$ are the spectrum and the
eigenvectors of the following generalized spectral problem:
\begin{equation} \label{m_eqn}
C^T \dot{f}_k + i \lambda_k C^T f_k = 0.
\end{equation}
\begin{proof}
For some $k\in \mathbb{N}$, we take $f_k \in H_0^1 \big( (0,T);Y
\big)$ such that $W^T f_k = u^{f_k} (T) = \varphi_k$. Since $f_k
(0) = f_k (T) = 0$ from our assumptions, the equalities
\[ u^{\frac{d}{dt} f_k} = \dfrac{d}{dt} u^{f_k}, \quad B f_k (T) = 0 \]
hold true. Then for arbitrary $g$, we can evaluate
\begin{align*}
\left( C^T \dfrac{d}{dt} f_k, g \right)_{\mathcal{F}^T} &= \left( u^{\frac{d}{dt} f_k} (T), v^g (T) \right)_H \\
&= \left( u_t^{f_k} (T), v^g (T) \right)_H\\
&= -i \left( Au^{f_k} (T) + B f_k (T), v^g (T) \right)_H \\
&= -i \left( A \varphi_k, v^g (T) \right)_H\\
&= -i (\lambda_k \varphi_k, v^g (T))_H\\
&= -i (\lambda_k u^{f_k} (T), v^g (T))_H\\
&= -i (\lambda_k C^T f_k, g)_{\mathcal{F}^T}.
\end{align*}
So the pairs $\{(\lambda_k, f_k)\}$ are solutions to
(\ref{m_eqn}). On the other hand, suppose that the pair $(\lambda,
f)$ is a solution to (\ref{m_eqn}) and $f \neq f_k$, $\lambda \neq
\lambda_k$ for all $k$. Then $W^T f$ has the form
\[ W^T f = u^f (T) = \sum_{k=1}^\infty a_k \varphi_k, \quad a_k \in \mathbb{R}.\]
We evaluate
\begin{align*}
0 = \left( C^T \dfrac{d}{dt} f_k + i \lambda C^T f, g \right)_{\mathcal{F}^T} &= \left( -i Au^f (T) + i \lambda W^T f, W_\#^T g \right)_H\\
&= -i \left( A \sum_{k=1}^\infty a_k \varphi_k - \lambda \sum_{k=1}^\infty a_k \varphi_k, W_\#^T g \right)_H \\
&= -i \left( \sum_{k=1}^\infty a_k (\lambda_k - \lambda) \varphi_k, W_\#^T g \right)_H.
\end{align*}
Using the spectral controllability assumption, we take $g = g_l$
such that $W_\#^T g_l = \psi_l$ for each $l$. Plugging this in the
right hand side of the above equality yields $a_k = 0$ for all $k$
and hence we obtain a contradiction. As a result, we have proved
the theorem.
\end{proof}
\end{theorem}

Similarly, one can find the set of controls for System (\ref{ED1}):
\begin{remark}
The spectrum of $A^*$ and (non-normalized) controls $g_k$ are the
spectrum and the eigenvectors of the following generalized
spectral problem:
\begin{equation} \label{m_eqn1}
\left(C^T \right)^* \dot{g}_k - i \overline{\lambda_k} \left( C^T \right)^* g_k = 0.
\end{equation}
\end{remark}


\section{Recovery of the Spectral Data in the Simple Case}

In this section, we will recover spectral data for System
(\ref{ED}). Let $u^f$ be the solution to System (2.1) and using
the Fourier method, we represent $u^f$ in the form
\[ u^f (t) = \sum_{k=1}^\infty c_k (t) \varphi_k, \quad c_k (t) = \int_0^t e^{-i \lambda_k (t-s)} \left( f(s), O \psi_k \right)_Y \: ds. \]
The response operator of the system is then given by
\[ \left( R^T f \right) (t) = \sum_{k=1}^\infty O \phi_k \int_0^t e^{-i \lambda_k (t-s)} \left( f(s), O \psi_k \right) \: ds. \]
These formulas motivate the following.
\begin{definition}
If $A$ satisfies Assumption \ref{A1}, then the set
\[ D := \left\{ \lambda_k, O \varphi_k, O \psi_k \right\}_{k=1}^\infty \]
is called the \emph{spectral data} of $A$.
\end{definition}

Having found eigenvalues $\{ \lambda_k \}_{k=1}^\infty$ and the sets of controls
$\{f_k\}_{k=1}^\infty$, $\{ g_k\}_{k=1}^\infty$ from Equations
(\ref{m_eqn}) and (\ref{m_eqn1}), we normalize the controls
according to the rule
\begin{equation}
\label{ContNorm} \left( C^T f_k, g_k \right)_{\mathcal{F}^T} = 1.
\end{equation}
Then, for $f_k$ and arbitrary $g$,
\begin{align} \label{TrEv}
\left( C^T f_k, g \right)_{\mathcal{F}^T} &= \left( \left( W_\#^T \right)^* W^T f_k,g\right)_{\mathcal{F}^T}\\
&= \left( \left(W_\#^T \right)^* \varphi_k, g\right)_{\mathcal{F}^T}. \notag
\end{align}
Our goal will be to evaluate $\left( W_\#^T \right)^* \varphi_k$
from the right hand side of (\ref{TrEv}).

Taking $a \in H$ we consider the system
\begin{equation} \label{ED4}
\left\{ \begin{array}{ll}
iw_t + Aw = 0, &\quad 0 < t < T,\\
w(T) = a,
\end{array} \right. \end{equation}
whose solution is denoted by $w^a$. We introduce the observation
operator for this system $\mathbb{O}^T$ by the rule
\begin{align*}
\mathbb{O}^T : \: &H \to L^2 \left( (0,T); Y \right)\\
&\left( \mathbb{O}^T a \right) (t) := \left( O w^a \right) (t).
\end{align*}
In particular, we provide the following lemma.
\begin{lemma} The observation operator $\mathbb{O}^T$ and $\left( W_\#^T \right)^*$ are related by
\begin{equation} \label{Wst}
\left( W_\#^T \right)^* = -i \mathbb{O}^T.
\end{equation}
\begin{proof}
Let $v^g$ be a solution to System (\ref{ED1}) and evaluate
\begin{align*}
\int_0^T \left( v_t^g (t), w^a (t) \right)_H \: dt &= \int_0^T i \left( A^* v^g + Bg, w^a \right)_H \: dt\\
&= \int_0^T i \left( \left( v^g, Aw^a \right)_H + \left( g, Ow^a\right)_Y \right) \: dt.
\end{align*}
On the other hand,
\begin{align*}
\int_0^T \left( v_t^g (t), w^a (t) \right)_H \: dt &= - \int_0^T \left( v^g, w_t^a \right)_H \: dt + \left( v^g, w^a \right)_H \big|_{t=0}^{t = T} \\
&= \int_0^T i \left( v^g, Aw^a \right)_H \: dt + \left( W_\#^T g,a \right)_H.
\end{align*}
Comparing the right hand sides of both equations yields
\[ i\int_0^T \left( g, Ow^a \right)_Y \: dt= \left( W_\#^T g, a \right)_H. \]
Since $g$ and $a$ were both arbitrary, the last equality completes
the proof of the lemma.
\end{proof}
\end{lemma}
We can infer from (\ref{Wst}) that $\left( W_\#^T \right)^*
\varphi_k = -i \mathbb{O}^T \varphi_k$ and setting $a = \varphi_k$
for System (\ref{ED4}) yields the solution $w^{\varphi_k} (t) =
\varphi_k e^{-i\lambda_k (T-t)}$. Hence, $\mathbb{O}^T \varphi_k =
O \varphi_k e^{-\lambda_k (T-t)}$. Plugging this into (\ref{TrEv})
gives that
\[ \left( C^T f_k, g \right)_{L^2 ((0,T);Y)}  = -i \left( O \varphi_k e^{-i \lambda_k (T-t)}, g \right)_{L^2 ((0,T);Y)}, \]
and thus
\begin{equation}
\label{TrPhi} -i O \varphi_k = e^{-i \lambda_k (t-T)} \left( C^T
f_k \right)(t) = \left( C^T f_k \right) (T).
\end{equation}
Similarly, it can be shown that
\begin{align*}
\left( f, \left( C^T \right)^* g_k \right)_{L^2 ((0,T);Y)} &= \left( f, \left( W^T \right)^* \psi_k \right)_{L^2 ((0,T);Y)}\\
&= \left( f, i O \psi_k e^{i \overline{\lambda_k} (T-t)} \right)_{L^2 ((0,T);Y)},
\end{align*}
and thus
\begin{equation}
\label{TrPsi} i O \psi_k = e^{i \overline{\lambda_k} (t-T)} \left(
\left( C^T \right)^* g_k \right) (t) = \left( \left( C^T \right)^*
g_k \right) (T).
\end{equation}
Hence, we propose the following method to calculate spectral data for
System (\ref{ED}) under Assumption \ref{A1}:

\begin{algorithm}\mbox{}
\begin{itemize}

\item[(1)] Solve generalized spectral problems (\ref{m_eqn}) and
(\ref{m_eqn1}) to find spectrum $\{\lambda_k\}_{k=1}^\infty$ and
controls $f_k,\,g_k$, $k=1,\ldots$.

\item[(2)] Normalize controls by (\ref{ContNorm}).

\item[(3)] Recover traces of eigenfunctions by (\ref{TrPhi}),
(\ref{TrPsi}).

\end{itemize}
\end{algorithm}


\section{The Spectral Problem and Recovery of the Spectral Data in the General Case}

We now assume that the operator $A$ satisfies the following:
\begin{assumption} \label{A2} \mbox{}
\begin{enumerate}[(a)]
\item The spectrum of $A$, denoted $\{ \lambda_k\}_{k=1}^\infty$,
is not simple. We denote the multiplicity of $\lambda_k$ by $L_k$.
\item The set of root vectors of $A$, $\{\varphi_k^l\}$, $k \in
\mathbb{N}$, $1 \leq l \leq L_k$, forms a Riesz basis in $H$. In
particular, for each $k \in \mathbb{N}$, the vectors in the chain
$\{\varphi_k^l\}_{l=1}^{L_k}$ satisfy
\begin{align*}
(A - \lambda_k I) \varphi_k^1 &= 0,\\
(A - \lambda_k I) \varphi_k^l &= \varphi_k^{l-1}, \quad 2 \leq l \leq L_k.
\end{align*}
\item The spectrum of $A^*$ is $\{
\overline{\lambda_k}\}_{k=1}^\infty$ and the root vectors of
$A^*$, $\{\psi_k^l\}$, $k \in \mathbb{N}$, $1 \leq l \leq L_k$,
also form a Riesz basis in $H$ and satisfy
\begin{align*}
(A - \overline{\lambda_k} I)\psi_k^{L_k} &= 0, \\
(A - \overline{\lambda_k} I) \psi_k^l &= \psi_k^{l+1}, \quad 1 \leq l \leq L_k -1.
\end{align*}
\item The property that $\left( \varphi_k^l, \psi_r^s \right)_H =
\delta_{kr} \delta_{ls}$ holds. \item Systems (\ref{ED}) and
(\ref{ED1}) are spectrally controllable. We denote the controls by
$f_k^l$ and $g_k^l$, both from $H_0^1 \big( (0,T);Y\big)$, such
that $W^T f_k^l = \varphi_k^l$ and $W_\#^T g_k^l = \psi_k^l$.
\end{enumerate}
\end{assumption}

The goal of this section is to obtain a result similar to Theorem
\ref{thm1}. In particular, we will construct generalized spectral problems from 
the spectra of $A$ and $A^*$ and the controls $\{f_k^l\}$, $\{g_k^l\}$
from Assumption \ref{A2}(e). We will also show that from these problems, 
we can obtain the spectra of $A$ and $A^*$ and normalized controls. 
We begin with the following lemma.
\begin{lemma} If $A$ satisfies Assumption \ref{A2}, then the spectrum of
$A$ and (non-normalized) controls $\{f_k^l$\} are solutions of the
following generalized spectral problem:
\begin{align} \label{CT_eqn}
C^T \dfrac{d}{dt} f_k^1 + i \lambda_k C^T f_k^1 &= 0,\\
C^T \dfrac{d}{dt} f_k^l + i \lambda_k C^T f_k^l &= -i C^T f_k^{l-1}, \quad 2 \leq l \leq L_k. \notag
\end{align}
\end{lemma}

\begin{proof}
Observe that
\begin{align*}
\left( C^T \dfrac{d}{dt} f_k^1, g \right)_{\mathcal{F}^T} &= \left( W^T \dfrac{d}{dt} f_k^1, W_\#^T g \right)_{\mathcal{F}^T}\\
&= \left( u_t^{f_k^1} (T), v^g (T) \right)_H\\
&= -i \left( A u^{f_k^1} (T) + Bf_k^1 (T), v^g (T) \right)_H\\
&= -i \left( A \varphi_k^1, v^g (T) \right)_H\\
&= -i \left( \lambda_k u^{f_k^1}(T), v^g (T) \right)_H\\
&= -i \left( \lambda_k C^T f_k^1, g \right)_{\mathcal{F}^T}.
\end{align*}
We then obtain the first equation of (\ref{CT_eqn}). We now let $2 \leq l \leq L_k$ and evaluate
\begin{align*}
\left( C^T \dfrac{d}{dt} f_k^l , g \right)_{\mathcal{F}^T} &= -i \left( Au^{f_k^l} (T), v^g (T) \right)_H\\
&= -i \left( \lambda_k \varphi_k^l + \varphi_k^{l-1}, v^g (T) \right)_H\\
&= -i \left( \lambda_k C^T f_k^l + C^T f_k^{l-1}, g \right)_{\mathcal{F}^T}.
\end{align*}
This yields the second equation. Hence, for each $k \in
\mathbb{N}$, the chain $\{f_k^l\}$ with eigenvalue 
$\lambda_k$ is a solution to the system of equations (\ref{CT_eqn}).
\end{proof}

Analogously, we have the following equations for the adjoint system (\ref{ED1}):
\begin{align} \label{CT*_eqn}
&\left(C^T\right)^* \dfrac{d}{dt} g_k^{L_k} - i \overline{\lambda_k} \left(C^T\right)^* g_k^{L_k} =0,\\
&\left(C^T\right)^* \dfrac{d}{dt} g_k^l - i \overline{\lambda_k} \left(C^T \right)^* g_k^l = i\left( C^T \right)^* g_k^{l+1}, \quad 1 \leq l \leq L_k -1. \notag
\end{align}

We will now consider the other direction. Let $\lambda$ be an
eigenvalue and $f^{(1)}, \ldots, f^{(M)}$ be a corresponding chain
of functions that solve Equation (\ref{CT_eqn}). For $f^{(1)}$, let
\[ W^T f^{(1)} = \sum_{k=1}^\infty \sum_{l=1}^{L_k} a_k^l \varphi_k^l, \quad a_k^l \in \mathbb{R}.\]
We now observe
\begin{align*}
0 &= \left( C^T \dfrac{d}{dt} f^{(1)} + i \lambda C^T f^{(1)}, g \right)_{\mathcal{F}^T} \\
&= -i \left( Au^{f^{(1)}} (T) - \lambda W^T f^{(1)}, W_\#^T g \right)_H\\
&= -i \left( \sum_{k=1}^\infty \sum_{l = 1}^{L_k} a_k^l (\lambda_k \varphi_k^l+ \varphi_k^{l-1}) - \lambda \sum_{k=1}^\infty \sum_{l=1}^{L_k} a_k^l \varphi_k^l, W_\#^T g \right)_H\\
&= -i \left( \sum_{k=1}^\infty \left[ \varphi_k^{L_k} a_k^{L_k} (\lambda_k - \lambda) + \sum_{l=1}^{L_k-1} \varphi_k^l \left( a_k^l (\lambda_k - \lambda_k)+ a_k^{l+1}\right) \right], W_\#^T g \right)_H.
\end{align*}
We again use spectral controllability of the system and see that
for every $k \in \mathbb{N}$ where $\lambda \neq \lambda_k$,
$a_k^l = 0$ for $1 \leq l \leq L_k$. In the case where $\lambda =
\lambda_k$, then $a_k^l = 0$ for $2 \leq l \leq L_k$ and $a_k^1$
is arbitrary. So $f^{(1)}$ is a (non-normalized) control that
drives System (\ref{ED}) to the state $\varphi_k^1$. We will
denote $a_k^1$ by $\alpha_1$ and express $W^T f^{(1)}$ as
\[ W^T f^{(1)} = \alpha_1 \varphi_k^1. \]

For $2 \leq j \leq M$, we will proceed by induction. Let
\begin{align*}
W^T f^{(j)} &= \sum_{r=1}^\infty \sum_{l=1}^{L_r} a_r^l \varphi_r^l,\\
W^T f^{(j-1)} &= \sum_{l=1}^{j-1} \alpha_{j-l} \varphi_k^l,
\end{align*}
and let $k \in \mathbb{N}$ be such that $\lambda = \lambda_k$. We evaluate
\begin{align*}
0 &= \left( C^T \dfrac{d}{dt} f^{(j)} + i \lambda_k C^T f^{(j)} + i C^T f^{(j-1)}, g\right)_{\mathcal{F}^T}\\
&= -i \bigg( \sum_{r=1}^\infty \sum_{j=1}^{L_r} a_r^j (\lambda_r \varphi_r^j + \varphi_r^{j-1}) - \lambda_k \sum_{r=1}^\infty \sum_{j=1}^{L_r} a_r^j \varphi_r^j - \sum_{l=1}^{j-1} \alpha_{j-l} \varphi_k^l, W_\#^T g \bigg)_H.
\end{align*}
For each $r \in \mathbb{N}$ with $\lambda_k \neq \lambda_r$, 
the corresponding summand has the form
\[ \varphi_r^{L_r} a_r^{L_r} (\lambda_r - \lambda_k) + \sum_{l=1}^{L_r -1} \varphi_r^l (a_r^l (\lambda_r - \lambda_k) + a_r^{l+1}). \]
Using the spectral controllability of the adjoint system, i.e., choosing 
$g = g_r^l$ with $W^T_\# g_r^l = \psi_r^l$, we
conclude that $a_r^l = 0$ for $1 \leq l \leq L_k$. However, 
for the case $\lambda_r=\lambda_k$, the corresponding summand has the form
\begin{align*}
\sum_{l=j}^{L_k-1} \varphi_k^l a_k^{l+1} + \sum_{l=1}^{j-1} \varphi_k^l \big( a_k^{l+1} - \alpha_{j-l} \big).
\end{align*}
Using spectral controllability yields the following results:
\begin{itemize}
\item $a_k^1$ is arbitrary and we will denote it by $\alpha_j$,
\item $a_k^l = \alpha_{j-l+1}$ for $2 \leq l \leq j$,
\item $a_k^l = 0$ for $j+1 \leq l \leq L_k$.
\end{itemize}
Hence, $W^T f^{(j)}$ has the form
\[ W^T f^{(j)} = \sum_{l=1}^j \alpha_{(j-l+1)} \varphi_k^l. \]

So, we obtain solutions to the generalized spectral problem,
$\{f_k^l\}$, $k \in \mathbb{N}$, $1 \leq l \leq L_k$, that are
controls to linear combinations of the root vectors associated to
$\lambda_k$. Similarly, for the generalized spectral
problem associated to the adjoint system, we obtain controls $\{g_k^l\}$ such that
\[ W_\#^T g_k^l = \sum_{j=1}^{L_k - l + 1} \beta_j \psi_k^{j+ l-1}. \]

As an example, suppose that $L_k = 3$, then the controls $\{f_k^l\}$ and $\{g_k^l\}$ are such that\\
\begin{minipage}{0.45\textwidth}
\begin{align*}
W^T f_k^1 &= \alpha_1 \varphi_k^1, \\
W^T f_k^2 &= \alpha_1 \varphi_k^2 + \alpha_2 \varphi_k^1,\\
W^T f_k^3 &= \alpha_1 \varphi_k^3 + \alpha_2 \varphi_k^2 + \alpha_3 \varphi_k^1.
\end{align*}
\end{minipage}
\begin{minipage}{0.45\textwidth}
\begin{align*}
W_\#^T g_k^1 &= \beta_1 \psi_k^1 + \beta_2 \psi_k^2 + \beta_3 \psi_k^3, \\
W_\#^T g_k^2 &= \beta_1 \psi_k^2 + \beta_2 \psi_k^3,\\
W_\#^T g_k^3 &= \beta_1 \psi_k^3.
\end{align*}
\end{minipage}

Our goal now is to construct a new family of controls to root vectors of $A$ and $A^*$ that are
biorthogonal in $H$ from the controls $\{f_k^l\}$, $\{g_k^l\}$. We first investigate the properties
of the control functions.
\begin{lemma} \label{ctrl1}
Let $2 \leq l,j \leq L_k$. Then $\big( C^T f_k^l, g_k^j \big) = \big( C^T f_k^{l-1}, g_k^{j-1} \big)$.
\begin{proof}
Evaluate
\begin{align*}
-i \big( C^T f_k^l, g_k^j \big)_{\mathcal{F}^T} &= \big( f_k^l, i (C^T)^* g_k^j \big)_{\mathcal{F}^T}\\
&= \bigg( f_k, (C^T)^* \dfrac{d}{dt} g_k^{j-1} - i \overline{\lambda_k} (C^T)^* g_k^{j-1} \bigg)_{\mathcal{F}^T}\\
&= \bigg( C^T \dfrac{d}{dt} f_k^l + i \lambda_k C^T f_k^l, g_k^{j-1} \bigg)_{\mathcal{F}^T}\\
&= -i \big( C^T f_k^{l-1}, g_k^{j-1} \big)_{\mathcal{F}^T}.
\end{align*}
Here, we use the assumption that the controls are from $H^1_0 ((0,T);Y)$.
\end{proof}
\end{lemma}
This lemma has two main applications. First, it shows the
connection between the controls in $\{f_k^l\}$, in particular, it
provides a strategy to inductively construct a biorthogonal
family. The second is
\begin{lemma} \label{ctrl2}
Let $2 \leq l,j \leq L_k$ with $l < j$. Then $\big(C^T f_k^l, g_k^j \big) = 0$.
\begin{proof}
Using the previous lemma, we observe
\begin{align*}
-i \big( C^T f_k^l, g_k^j \big)_{\mathcal{F}^T} &= -i \big(C^T f_k^1, g_k^{j-l+1} \big)_{\mathcal{F}^T}\\
&= \bigg( C^T \dfrac{d}{dt} f_k^1 + i \lambda_k C^T f_k^1, g_k^{j-l} \bigg)_{\mathcal{F}^T}\\
&= 0.
\end{align*}
\end{proof}
\end{lemma}

We will demonstrate how to construct controls that produce a
biorthogonal family. Suppose that we have solved (\ref{CT*_eqn})
and obtained $\{g_k^j\}$ for some $k \in \mathbb{N}$ and all $1
\leq j \leq L_k$. Let $f_k^1$ be a solution to the first equation
in (\ref{CT_eqn}). We then obtain $\hat{f}_k^1$ by the rule
\begin{equation} \label{fhat1} \hat{f}_k^1 = \dfrac{1}{ \big( C^T f_k^1, g_k^1 \big)_{\mathcal{F}^T}} f_k^1.  \end{equation}
We note that $\big(C^T \hat{f}_k^1, g_k^1 \big)_{\mathcal{F}^T} =
1$ and for $2 \leq j \leq L_k$, $\big( C^T \hat{f}_k^1, g_k^j
\big)_{\mathcal{F}^T} = 0$ as a result of Lemmas \ref{ctrl1} and \ref{ctrl2}. We then
use $\hat{f}_k^1$ to solve (\ref{CT_eqn}) and obtain $f_k^2$. We 
construct $\hat{f}_k^2$ by the rule
\[ \hat{f}_k^2 = f_k^2 - \big( C^T f_k^2, g_k^1 \big)_{\mathcal{F}^T} \hat{f}_k^1. \]
By construction, $\big(C^T \hat{f}_k^2, g_k^j
\big)_{\mathcal{F}^T} = \delta_{2,j}$ for $2 \leq j \leq L_k$ and
$\big(C^T \hat{f}_k^2, g_k^1\big)_{\mathcal{F}^T} = 0$.

We will now proceed iteratively. Let $f_k^l$ be the solution
obtained from $\hat{f}_k^{l-1}$ and define $\hat{f}_k^l$ by
\begin{equation} \label{fhatl}
\hat{f}_k^l = f_k^l - \big(C^T f_k^l, g_k^1 \big)_{\mathcal{F}^T} \hat{f}_k^1.
\end{equation}
In this way, we obtain a new collection $\{\hat{f}_k^l\}$ such that
\[ \big( C^T \hat{f}_k^l, g_k^j \big)_{\mathcal{F}^T} = \delta_{lj}. \]
We define two new collections of root vectors $\{
\hat{\varphi}_k^l\}$, $\{\hat{\psi}_k^l\}$, for $k \in
\mathbb{N}$, $1 \leq l \leq L_k$ where
\[ \hat{\varphi}_k^l = W^T \hat{f}_k^l, \qquad
\hat{\psi}_k^l = W_\#^T g_k^l. \] We note that this is a
biorthogonal family and each collection is also a Riesz basis in
$H$. To calculate the spectral data of $A$ under Assumption \ref{A2}, we propose the 
following method:
\begin{algorithm} \label{alg2} \mbox{}
\begin{itemize}

\item[(1)] Solve the generalized spectral problem
(\ref{CT*_eqn}) to find the spectrum $\{\lambda_k\}_{k=1}^\infty$ and
controls $\{g_k^j\}$, for $k\in \mathbb{N}$.

\item[(2)] Solve the generalized spectral problem
(\ref{CT_eqn}) for $f_k^1$ and construct $\hat{f}_k^1$ 
according to (\ref{fhat1}).

\item[(3)] Iteratively solve (\ref{CT_eqn}) using $\hat{f}_k^{l-1}$ to obtain
$f_k^l$ and construct $\hat{f}_k^l$ by (\ref{fhatl}).

\item[(4)] Recover traces of eigenfunctions by (\ref{TrPhi}),
(\ref{TrPsi}).

\end{itemize}
\end{algorithm}


\section{Dynamical Inverse Problem for the nonsymmetric matrix Schr\"odinger Operator on an Interval}

We now consider the problem of recovering a (nonsymmetric) matrix
potential $Q$ of the following dynamical system
\begin{equation}
\label{Dyn_syst1}
\left\{ \begin{array}{ll}
iu_t - u_{xx} + Q(x)u = 0, &\quad 0 \leq x \leq \ell, \: 0 < t < T,\\
u(x,0) =  0, &\quad 0 \leq x \leq \ell,\\
u(0,t) = f(t), \: u(\ell,t) = 0, &\quad 0 < t < T,
\end{array} \right.
\end{equation}
from the \emph{response operator} defined by
$(R^Tf)(t):=u^f_x(0,t)$. Recall that $u$ is a vector-valued function.
We denote the space of controls by $\mathcal{F}^T = L^2
\big((0,T);\mathbb{R}^N\big)$.  We will first show that the
Schr\"odinger equation is null controllable, which is equivalent
to exact controllability, and hence spectrally controllable. We
will prove this using the control transmutation method (see
\cite{miller2005}, Sections 8 and 9). Afterwards, we will use our previous
results to recover the spectral data of the system. We then use
the spectral data to recover the matrix potential $Q$.

\subsection{Spectral Controllability of the Schr\"odinger Equation with nonsymmetric matrix potential}\mbox{}

Consider the following system:
\begin{equation} \label{nullSchr}
\left\{ \begin{array}{ll}
iu_t - u_{xx} + Q(x)u = 0, &\quad 0 \leq x \leq \ell, \: 0 < t < T,\\
u(x,0) = \varphi_0 (x),&\quad 0 \leq x \leq \ell,\\
u(0,t) = f(t), \: u(\ell,t) = 0, &\quad 0 < t < T,
\end{array} \right.
\end{equation}
where $\varphi_0 \in L^2 ((0,\ell);\mathbb{R}^N)$ is the initial state of the system. 
We will show null controllability of this system, i.e., the existence of 
a control function $f$ such that the corresponding solution
satisfies $u(x,T) = 0$. We begin by constructing the auxiliary
wave system
\begin{equation} \label{nullWave}
\left\{ \begin{array}{ll}
v_{tt} - v_{xx} + Q(x)v = 0, &\quad 0 \leq x \leq \ell, \: 0 < t < T^*,\\
v(x,0) = \varphi_0 (x), \: v_t (x,0) = 0, &\quad 0 \leq x \leq \ell,\\
v(0,t) = g(t), \: v(\ell,t) = 0, &\quad 0 < t < T^*.
\end{array} \right.
\end{equation}
It is a known result (see \cite{ab1996}) that for $T^* \geq 2\ell$, 
System (\ref{nullWave}) is exactly controllable and
hence null controllable.

Let $\varphi_0$, $L> 0$, and $T>0$ be given. We choose $T^* \geq 2\ell$ 
and thus System (\ref{nullWave}) is null
controllable. Hence, for System (\ref{nullWave}), we obtain the
control function $g(t)$ and the solution $v(x,t)$ with the
property that
\[ v(x,T^*) = v_t (x,T^*) = 0. \]
We now extend $v$ and $g$ to $\tilde{v}$ and $\tilde{g}$ by the rule
\begin{align*}
\tilde{v} (x,-t) &= \tilde{v}(x,t) = v(x,t), \\
\tilde{g} (-t) &= \tilde{g} (t) = g(t),
\end{align*}
for $0 < t < T^*$. We note that $\tilde{v}$ inherits the following
properties from $v$:
\begin{equation} \label{tildevxt} \tilde{v} (x,-T^*) = \tilde{v}(x,T^*) = \tilde{v}_t (x,-T^*) = \tilde{v}_t (x,T^*) = 0. \end{equation}

We define a scalar function $k(s,t)$ to be the solution to the system
\begin{equation} \left\{ \begin{array}{ll} \label{kst}
i \partial_t k - \partial_s^2 k = 0, &\quad -T^* \leq s \leq T^*, \: 0 \leq t \leq T,\\
k (s,0) = \delta (s), \: k(s,T) = 0.
\end{array} \right. \end{equation}
The existence of $k(s,t)$ is a result of System (\ref{kst}) being exactly
controllable from both ends (see \cite{miller2005} Section 2), we omit the
boundary conditions as they do not need to be specified for our
purposes. We need only its initial state and its state at time $t
= T$. We then construct $f$ and $u$ by
\begin{align*}
f(t) &= \int_{-T^*}^{T^*} k(s,t) \tilde{g} (s) \: ds,\\
u(x,t) &= \int_{-T^*}^{T^*} k(s,t) \tilde{v} (x,s) \: ds.
\end{align*}
We observe that $u$ inherits the following properties:
\begin{align*}
u(x,0) &= \varphi_0 (x), \\
u(0,t) &= f(t), \\
u(L, t) &= 0,\\
u(x,T) &= 0.
\end{align*}
These properties, along with (\ref{tildevxt}), demonstrates that
$u$ is a solution to System (\ref{nullSchr}) with control
$f$, and this proves that the system is null controllable.

\subsection{Recovery of the Spectral Data and the Matrix Potential}\mbox{}

Returning to System (\ref{Dyn_syst1}), we construct the connecting
operator, $C^T$, from $R^T$ by means of (\ref{CT_repr}). We then
implement Algorithm \ref{alg2} to obtain the eigenvalues
$\{\lambda_k\}$, and controls $\{\hat{f}_k^l\}$, $\{g_k^l\}$
for $k \in \mathbb{N}$ and $1 \leq l \leq L_k$, where $L_k$ is the
multiplicity of the eigenvalue $\lambda_k$. Note that
\begin{align*}
R^T \hat{f}_k^l &= \dfrac{d}{dx} \hat{\varphi}_k^l (x) \bigg|_{x = 0} =: \Phi_{k,l},\\
R^T_\# g_k^l &= \dfrac{d}{dx} \hat{\psi}_k^l (x) \bigg|_{x=0} =: \Psi_{k,l}.
\end{align*}
So we have obtained spectral data, $\{\lambda_k, \Phi_{k,l}, \Psi_{k,l}\}$, for System (\ref{Dyn_syst1}). \\

To recover the matrix potential, we construct the auxiliary wave system
\begin{equation}
\label{Dyn_syst} \left\{
\begin{array}{ll}
w_{tt}-w_{xx}+Q(x)w=0, &\quad 0 \leq x\leq\ell, \: 0<t<T,\\
w(x,0)=w_t(x,0)=0,&\quad 0\leq x\leq \ell,\\
w(0,t)=f(t),\: w(\ell,t)=0, &\quad 0<t<T,
\end{array}
\right.
\end{equation}
with response operator $(R^T_w f) (t) := w_x^f (0,t)$. Our next
step is to express the connecting operator, $C^T_w$, for System
(\ref{Dyn_syst}) in terms of the spectral data obtained from
(\ref{Dyn_syst1}). Using the Fourier method, we represent the
solution, $w(x,t)$, in the form
\[ w^f (x,t) = \sum_{k,l} b_k^l (t) \hat{\varphi}_k^l (x), \]
where
\begin{align} \label{bkl_def}
b_k^1 (t)&= \int_0^t \Big[ \Psi_{k,1} f(\tau) \Big] \dfrac{\sin \sqrt{\lambda_k} (t-\tau)}{\sqrt{\lambda_k}},\\
 b_k^l (t) &= \int_0^t \Big[ \Psi_{k,l} f(\tau) - b_k^{l-1} (\tau) \Big] \dfrac{\sin \sqrt{\lambda_k} (t-\tau)}{\sqrt{\lambda_k}} \: d\tau, \quad 2 \leq l \leq L_k. \notag
\end{align}
Similarly, we denote $w_\#^g (x,t)$ to be the solution to the
adjoint wave system with the representation
\[ w_\#^g (x,t) = \sum_{k,l} c_k^l (t) \psi_k^l (x), \]
where
\begin{align} \label{ckl_def}
c_k^1 (t) &= \int_0^t \Big[ \Phi_{k,1} g(s) \Big] \dfrac{\sin \sqrt{\overline{\lambda_k}} (t-s)}{\sqrt{\overline{\lambda_k}}},\\
c_k^l (t) &= \int_0^t \Big[ \Phi_{k,l} g(s) - c_k^{l-1} (s) \Big] \dfrac{\sin \sqrt{\overline{\lambda_k}} (t-s)}{\sqrt{\overline{\lambda_k}}}, \quad 2 \leq l \leq L_k. \notag
\end{align}
From here, we compute
\begin{align} \label{CTw_repr}
(C_w^T f, g)_{\mathcal{F}^T} &= \sum_{k,l} \Bigg\{ \int_0^T \Big[ \Psi_{k,l} f(t) - b_k^{l-1} (t) \Big] \dfrac{\sin \sqrt{\lambda_k} (T-t)}{\sqrt{\lambda_k}} \: dt \cdot \\
&\qquad \cdot \int_0^T \Big[ \Phi_{k,l} g(s) - c_k^{l-1} (s)\Big] \dfrac{\sin \sqrt{\overline{\lambda_k}}(T-s)}{\sqrt{\overline{\lambda_k}}} \: ds \Bigg\}, \notag
\end{align}
and define $b_k^0 (t) \equiv 0$, $c_k^0 (t) \equiv 0$. Hence, $C_w^T$ is completely determined by the spectral data. \\

Now let $y_j (x)$ be the solution to the boundary value problem
\begin{equation}\label{y_eqn}
\left\{ \begin{array}{l} y''(x) - Q(x) y(x) = 0, \quad 0 \leq x \leq \ell,\\
y(0) = 0, \: y'(0) = e_j,
\end{array} \right. \end{equation}
where $e_j$ is the $j$-th standard basis vector in $\mathbb{R}^N$.
Let $p_j^T$ be the control function such that
\begin{equation} \label{wpt_eqn}
w^{p_j^T} (x,T) = \left\{ \begin{array}{rl} y(x),& \quad x \leq T,\\
0, &\quad x > T. \end{array} \right.
\end{equation}
For any $g \in C_0^\infty ((0,T);\mathbb{C}^N)$, we have
\begin{align*}
(C_w^T p_j^T, g)_{\mathcal{F}^T} &= (w^{p_j^T} (\cdot, T), w_\#^g (\cdot, T))_{L^2 ((0,T);\mathbb{R}^N)}\\
&= \int_0^T \langle y_j(x), w_\#^g (x,T)\rangle \: dx\\
&= \int_0^T (T-t) \: dt \int_0^T \langle y_j(x), (w_\#^g)_{tt} (x,t) \rangle \: dx\\
&= \int_0^T (T-t) \: dt \int_0^T \bigg\langle y_j(x), \Big[ (w_\#^g)_{xx} (x,t) - Q^*(x)w_\#^g (x,t)\Big] \bigg\rangle \: dx\\
&= \int_0^T (T-t) \: dt \Bigg\{ \Big\langle y_j''(x) + Q(x) y_j(x), w_\#^g (x,t) \Big\rangle \: dx \\
&\qquad + \bigg[ \Big\langle y_j(x), (w_\#^g)_x (x,t)\Big\rangle  - \Big\langle y_j'(x), w_\#^g (x,t)\Big\rangle \bigg]_{x=0}^{x=T} \Bigg\}\\
&= \int_0^T \Big\langle (T-t)e_j, g(t) \Big\rangle \: dt.
\end{align*}
In the previous calculation, we use that for $g \in C_0^\infty
((0,T);\mathbb{C}^N)$, the function $w_\#^g$ and its derivatives
are equal to zero at $x = T$. Hence, the function $p_j^T$
satisfies the equation
\begin{equation*}
\big( C^T p_j^T \big) (t) = (T-t)e_j, \quad t \in (0,T).
\end{equation*}
Since $C_w^T$ is boundedly invertible, this equation has a unique
solution, $p_j^T \in \mathcal{F}^T$, for any $T \leq N$. Moreover,
it can be proved that $p_j^T \in H^1 ((0,T); \mathbb{C}^N)$ and
\begin{equation}
w^{p_j^T} (T-0,T) = -p_j^T (+0) =: - \mu_j (T),
\end{equation}
(see, for example, \cite{alp2002,alp2005}). From (\ref{wpt_eqn}),
$w^{p_j^T} (T-0,T) = y_j(T)$ and thus $\mu_j (T)$ is twice
differentiable with respect to $T$. We then construct the $N
\times N$ matrix $M (T)$ by
\begin{equation*}
M(T) = \left[ \begin{array}{c|c|c|c}
\mu_1 (T) & \mu_2 (T) & \cdots & \mu_n (T) \end{array} \right]
= \left[ \begin{array}{c|c|c|c} y_1 (T) & y_2 (T) & \cdots & y_n (T) \end{array} \right].
\end{equation*}
The matrix $M(T)$ is invertible except for finitely many points
(see \cite{abi1991}) and from (\ref{y_eqn}) we obtain
\[ Q(T)= M''(T) M^{-1}(T). \]
By varying $T$ in $(0,\ell)$, we obtain $Q(\cdot)$ in that
interval. For the finite number of times that $M(T)$ is singular,
we recover $Q(T)$ by continuity. This completes the process of
recovering the nonsymmetric potential matrix $Q$ on an interval.


\section*{Acknowledgments} Sergei Avdonin  was  supported in part by NSF grant DMS
1909869 and by the Ministry of Education and Science of Republic
of Kazakhstan under the grant No. AP05136197.  
Alexander Mikhaylov was  supported in part by in part RFBR 17-01-00099 and RFBR
18-01-00269.
Victor Mikhaylov
was supported in part by RFBR 17-01-00529, RFBR 18-01-00269.


\bibliographystyle{ieeetr}

\end{document}